\documentclass{amsart}
\usepackage{hyperref}
\newtheorem{thm}{Theorem}[section]
\newtheorem{lem}[thm]{Lemma}
\newtheorem{cor}[thm]{Corollary}
\newtheorem{conj}[thm]{Conjecture}
\newtheorem{prop}[thm]{Proposition}
\theoremstyle{definition}

\newtheorem{rem}[thm]{Remark}
\numberwithin{equation}{section}

\def\Q{{\mathbf Q}}
\def\Z{{\mathbf Z}}
\def\C{{\mathbf C}}
\def\R{{\mathbf R}}

\def\Gal{\mathrm{Gal}}
\def\End{\mathrm{End}}

\def\Hom{\mathrm{Hom}}
\def\Jac{\mathrm{Jac}}
\def\genus{\mathrm{genus}}

\def\cond{\mathrm{cond}}

\def\rank{\mathrm{rank}}

\def\eqref#1{(\ref{#1})}

\title[Ranks of quadratic twists of elliptic curves]
{Rank frequencies for quadratic twists of elliptic curves}

\thanks{We thank NSF, NSA, and the Alexander-von-Humboldt
Stiftung for financial support, and AIM and 
the Mathematics Institute of the University of Erlangen 
for congenial working environments.}

\author[K.\ Rubin]{K.\ Rubin}
\address{Department of Mathematics, Stanford University, 
Stanford, CA 94305}
\email{rubin\char`\@math.stanford.edu}

\author[A.\ Silverberg]{A.\ Silverberg}
\address{Department of Mathematics, Ohio State University, 
Columbus, Ohio 43210}
\email{silver\char`\@math.ohio-state.edu}

\begin{document}

\begin{abstract}
We give explicit examples of infinite families of elliptic curves 
$E$ over ${\bf Q}$ with (nonconstant)
quadratic twists over ${\bf Q}(t)$ of rank at least 
$2$ and $3$. We recover some results announced by Mestre, as well as some
additional families.
Suppose $D$ is a squarefree integer and let $r_E(D)$ denote the rank of 
the quadratic twist of $E$ by $D$.
We apply results of Stewart and Top to our examples
to obtain results of the form  
$$\#\{D : |D|<x,  \,\,  r_E(D) \ge 2\} \gg x^{1/3}$$
$$\#\{D : |D|<x, \,\,  r_E(D) \ge 3\} \gg x^{1/6}$$
for all sufficiently large $x$.
\end{abstract}

\maketitle

\section{Introduction}
\label{intro}

Throughout this paper $E$ is an elliptic curve over $\Q$ defined by 
a Weierstrass equation $y^2 = f(x)$, where $f$ is a monic cubic
polynomial.  The curve $Dy^2 = f(x)$ will be denoted $E_D$.
When $D$ is a nonzero integer, let $r_E(D)$ denote the rank of $E_D(\Q)$. 
Let 
\begin{align*}
N_r(E,x) &= \#\{\text{squarefree $D\in\Z$} : 
    \text{$|D|<x$ and $r_E(D) \ge r$}\}, \\
N_{r}^{+}(E,x) &= \#\{\text{squarefree $D\in\Z$} : 
    |D|<x, r_E(D) \ge r, r_E(D)\equiv r \hskip-5pt\pmod{2}\}.
\end{align*}
In \cite{Gouvea-Mazur}, Gouv\^ea and Mazur showed (using the fact that 
the twist $E_{f(u)}$ has rank one over $\Q(u)$) 
that if the Parity Conjecture holds then 
$$N_{2}^+(E,x) > x^{\frac{1}{2}-\epsilon}$$
for all sufficiently large $x$.

In Th\'eor\`eme 1 of \cite{Mestre}  Mestre showed that 
if $j(E) \notin \{0,1728\}$ then 
there is a polynomial $g(u) \in \Q[u]$ of degree $14$ such that the twist
$E_{g(u)}$ has rank at least $2$ over $\Q(u)$.
In Theorem 3 of \cite{STop}, Stewart and Top used Mestre's result to
show that 
$$N_{2}(E,x) \gg x^{\frac{1}{7}}/(\log x)^2$$
for such $E$ and for all sufficiently large $x$. 
For a special family of elliptic curves 
$E$, using a twist of $E$ over $\Q(u)$ of rank at least $3$, 
Stewart and Top (Theorem 6 of \cite{STop}) 
found lower bounds for $N_{3}(E,x)$.
Mestre announced in Th\'eor\`eme 2 of \cite{Mestre98} 
that if the torsion subgroup
of $E(\Q)$ contains $\Z/2\Z \times \Z/4\Z$,
then $E$ has a (nonconstant) quadratic twist over $\Q(u)$ of rank at least $3$.

For certain elliptic curves $E$, Howe, Lepr\'evost, and Poonen 
(see Proposition 4 of \cite{HLP}) 
constructed polynomials $g(u)$ of degree $6$ such that the twist $E_{g(u)}$
has rank $2$ over $\Q(u)$. 

In this paper we describe a method (\S\ref{useful}) 
for constructing (nonconstant) quadratic twists of $E$ over $\Q(u)$
of ranks (at least) $2$ and $3$, and obtain further examples.  
In the rank $2$ case (\S\ref{rank2}) we show that this method 
recovers the above mentioned
results of Howe, Lepr\'evost, and Poonen and of Mestre. 
The rank $3$ cases (\S\ref{rank3}) include Mestre's curves and 
some other infinite families.  
In \S\ref{densect} we use results of Stewart and Top to obtain lower
bounds for $N_r(E,x)$ (and for $N_{r+1}^+(E,x)$, 
subject to the Parity Conjecture) for these examples, with $r=2$ or $3$.

The idea behind the method is that given an elliptic curve $E$ over $\Q(t)$, 
it is easy to find twists of $E$ of rank $r$ over extensions $K/\Q(t)$ 
with $\Gal(K/\Q(t)) \cong (\Z/2\Z)^{r-1}$.  When $r \le 3$, we show how to 
do this with $K = \Q(u)$ for some $u$, for certain families of curves.  

We used PARI and Mathematica to perform the computations in this paper.
The results of the computations, 
including those which are too long to display in the paper,
can be found in the electronic appendix \cite{ElecApp}.

We would like to thank Jean-Fran{\c c}ois Mestre for pointing out that
the curves with $(\Z/2\Z \times \Z/8\Z)$-torsion are isogenous to twists
of the curves in Theorem 6 of \cite{STop}, 
and Brian Conrey for telling us about connections between 
rank heuristics coming from Random Matrix Theory and
Theorems \ref{hlpcor} and \ref{hlpcor2} below.

After writing this paper we learned that the method we use here to 
construct rank $2$ and $3$ quadratic twists is essentially the same as 
one of the methods used by Mestre to prove the results 
announced in \cite{Mestre98}.
Since Mestre's proofs and explicit descriptions
of the twists he obtains have not been published, and we need  
explicit forms of these twists for the applications in \S\ref{densect}, 
we include the details here. 

\section{Constructing useful twists}
\label{useful}

We begin with the following well-known result.

\begin{lem}
\label{basic}
If $F$ is a field of characteristic different from $2$, 
$A$ is an elliptic curve over $F$, and 
$K$ is an abelian extension of $F$ with $\Gal(K/F) \cong (\Z/2\Z)^d$, then 
$$
\rank(A(K)) = \sum_{\chi}\rank(A^{\chi}(F))
$$
where the sum is over characters $\chi : \Gal(K/F) \to \{\pm 1\}$, 
and $A^{\chi}$ is $A$ if $\chi = 1$ and otherwise $A^{\chi}$ is 
the quadratic twist of $A$ corresponding to $\chi$.
\end{lem}

\begin{cor}
\label{method}
Suppose $E$ is an elliptic curve over $\Q$,  
$g_1,\ldots,g_r \in \Q(t)^\times$,
the fields $\Q(t,\sqrt{g_i})$ are distinct quadratic extensions
of $\Q(t)$, and
$\rank(E_{g_i}(\Q(t))) > 0$ for every $i$.
Then $$\rank(E_{g_1}(\Q(t,\sqrt{g_1g_2},\ldots,\sqrt{g_1g_r}))) \ge r.$$
If in addition $\Q(t,\sqrt{g_1g_2},\ldots,\sqrt{g_1g_r}) = \Q(u)$ for 
some $u$, and $g(u) = g_1(t)$, then $\rank(E_{g(u)}(\Q(u))) \ge r$.
\end{cor}

\begin{proof}
Take $A = E_{g_1}$, 
$F = \Q(t)$, and $K = \Q(t,\sqrt{g_1g_2},\ldots,\sqrt{g_1g_r})$.
By Lemma \ref{basic},
\begin{align*}
\rank(E_{g_1}(\Q(t,\sqrt{g_1g_2},\ldots,&\sqrt{g_1g_r}))) \\
    &\ge \rank(E_{g_1}(\Q(t))) + \sum_{i=2}^r \rank(E_{g_1(g_1g_i)}(\Q(t))) \\
    &= \sum_{i=1}^r \rank(E_{g_i}(\Q(t))) \ge r.
\end{align*}
This proves the first part of the corollary, and the second is immediate.
\end{proof}

Given an elliptic curve $E$ over $\Q$, 
we want to use Corollary \ref{method} to construct twists of 
$E$ over $\Q(u)$ of ``large'' rank.  The following lemma provides us with 
elements $g \in \Q(t)$ such that $\rank(E_g(\Q(t))) > 0$.

\begin{lem}
\label{gm}
Suppose $E$ is the elliptic curve over $\Q$ defined by $y^2 = f(x)$.
Then for every nonconstant $h \in \Q(t)$ we have
$$
\rank(E_{f(h(t))}(\Q(t))) > 0.
$$
\end{lem}

\begin{proof}
The point $(h(t),1)$ belongs to $E_{f(h(t))}(\Q(t))$.  Since this 
point is nonconstant, it cannot be a torsion point. 
\end{proof}

\begin{rem}
Conversely, if $g \in \Q(t)$ and $\rank(E_g(\Q(t))) > 0$,
then there is an $h \in \Q(t)$ such that $E_g \cong E_{f(h(t))}$.  
To see this, let $(h,k)$ be a point of infinite order in $E_g(\Q(t))$, 
and then $f\circ h = k^2 g$.
\end{rem}

To apply Corollary \ref{method} we also need to know when 
$\Q(t,\sqrt{g_1g_2},\ldots,\sqrt{g_1g_d})$ is a 
rational function field.  For this we use the following well-known 
result.

\begin{lem}
\label{hypgenus}
Suppose $k \in \Q[t]$ is a nonconstant squarefree polynomial.  
Then the curve $s^2 = k(t)$ has genus $\bigl[\frac{\deg(k)-1}{2}\bigr]$.
\end{lem}

\begin{cor}
\label{genuscor}
\begin{enumerate}
\item
If $k \in \Q[t]$ is squarefree and $1 \le \deg(k) \le 2$, then 
the function field $\Q(t,\sqrt{k})$ has genus zero.
\item
If $k_1, k_2 \in \Q[t]$ are linear and linearly independent over $\Q$, then 
the function field $\Q(t,\sqrt{k_1},\sqrt{k_2})$ has genus zero.
\end{enumerate}
\end{cor}

\begin{proof}
The first statement is immediate from Lemma \ref{hypgenus}.  
The second statement follows without difficulty by applying 
(i) first to the extension $\Q(t,\sqrt{k_1})/\Q(t)$, and 
then to the extension $\Q(t,\sqrt{k_1},\sqrt{k_2})/\Q(t,\sqrt{k_1})$.
\end{proof}

If $g(t) \in \Q(t) \subseteq \Q(u)$, then $g(u) \in \Q(u)$ will denote the 
element $g(t(u))$, where $t(u)$ is the image of $t$ in $\Q(u)$.
We regard $f$ as an element of $\Q[t]$.

The following two propositions summarize a
method for producing twists of $E$ over $\Q(u)$ 
with ranks (at least) $2$ and $3$.

\begin{prop}
\label{rk2}
Suppose $h \in \Q(t)$ is such that $f\circ h = k f j^2$ 
with $j \in \Q(t)$, $k \in \Q[t]$, and $k$ squarefree. 
If $\deg(k)=1$, then the function 
field $\Q(t,\sqrt{k(t)}) = \Q(u)$ with $u=\sqrt{k(t)}$, 
and we have $\deg(f(u))=6$ and 
$\rank(E_{f(u)}(\Q(u))) \ge 2$.
If $\deg(k)=2$ and the curve $s^2 = k(t)$ has a rational point, 
then $\Q(t,\sqrt{k}) = \Q(u)$ for some $u$, and 
$\rank(E_{f(u)}(\Q(u))) \ge 2$.
\end{prop}

\begin{proof}
This follows directly from Corollary \ref{method} (with $g_1 = f$
and $g_2 = f\circ h$), Lemma \ref{gm}, and Corollary \ref{genuscor}.
\end{proof}

\begin{prop}
\label{rk3}
Suppose $h_1, h_2 \in \Q(t)$ are such that $f\circ h_i = k_i f j_i^2$ 
for $i = 1, 2$, with $j_i \in \Q(t)$, $k_i \in \Q[t]$, and $k_i$ linear and 
$\Q$-linearly independent.  
If the curve $s^2 = k_1(t), r^2 = k_2(t)$ has a rational point, 
then the function field $\Q(t,\sqrt{k_1},\sqrt{k_2}) = \Q(u)$ for some 
$u$, and $\rank(E_{f(u)}(\Q(u))) \ge 3$.
\end{prop}

\begin{proof}
This follows directly from Corollary \ref{method}, 
Lemma \ref{gm}, and Corollary \ref{genuscor}.
\end{proof}

To apply Propositions \ref{rk2} or \ref{rk3}, we want to find elements 
$h \in \Q(t)$ such that $f\circ h = k f j^2$ with $j \in \Q(t)$,
$k \in \Q[t]$, and 
$k$ linear.  The following two propositions give two possible ways of 
doing this.

\begin{prop}
\label{lfts}
Suppose $h(t) = \frac{\alpha t + \beta}{t + \delta} \in \Q(t)$ 
is a linear fractional transformation which permutes the roots of $f$.
Then 
$$
f(h(t)) = f(\alpha)(t + \delta) f(t) (t + \delta)^{-4}.
$$
\end{prop}

\begin{proof}
Both sides have the same divisor, and evaluate to 
$f(\alpha)$ at $t = \infty$.
\end{proof}

\begin{rem}
\label{isogsrem}
Suppose $\tilde{E}$ is an elliptic curve $Y^2 = \tilde{f}(X)$ with
$\tilde{f}$ a monic cubic, and suppose 
$\phi : \tilde{E} \to E$ is an isogeny.  Then $\phi(X,Y) = (\phi_x(X),Y\phi_y(X))$ 
with $\phi_x, \phi_y \in \Q(t)$, since the $x$-coordinate of 
$\phi$ is an even function on $\tilde{E}$ and the $y$-coordinate is an odd 
function.  
\end{rem}

\begin{prop}
\label{isogs}
Suppose $\tilde{E}$ is an elliptic curve $Y^2 = \tilde{f}(X)$ with
$\tilde{f}$ a monic cubic, and suppose
$\phi : \tilde{E} \to E$ is an isogeny.  
If $\phi_x$ and $\phi_y$ are as in Remark \ref{isogsrem},
$\mu(t) = \frac{\alpha t + \beta}{t + \delta} \in \Q(t)$ 
is a linear fractional transformation which sends the roots of $f$ to the 
roots of $g$, 
and $h(t) = \phi_x(\mu(t))$, then 
$$
f(h(t)) = 
\tilde{f}(\alpha)(t + \delta) f(t) 
\Bigl(\frac{\phi_y(\mu(t))}{(t + \delta)^{2}}\Bigr)^2.
$$
\end{prop}

\begin{proof}
By Remark \ref{isogsrem}, 
$
f(\phi_x(X)) = Y^2\phi_y(X)^2 = \tilde{f}(X)\phi_y(X)^2$.
As in the proof of Proposition \ref{lfts}, 
$$
\tilde{f}(\mu(t)) = \tilde{f}(\alpha)(t + \delta) f(t) (t + \delta)^{-4}
$$
and the identity of the proposition follows.
\end{proof}

\begin{rem}
\label{hlprem}
Suppose $g(u)\in\Q[u]$ is squarefree and
nonconstant, and let $C$ be the curve
$s^2=g(u)$. Then
$$\rank(E_g(\Q(u))) = \rank(\Hom_\Q(\Jac(C),E)) \le \genus(C)$$
(see \S 4 of \cite{STop}).
\end{rem}

\section{Rank 2}
\label{rank2}

The following statement is a reformulation of a result of
Howe, Lepr\'evost, and Poonen (Proposition 4 of \cite{HLP})
in a special case.  The proof below is different from theirs, 
and uses the method described in the preceding sections.

\begin{thm}[\cite{HLP} Proposition 4]
\label{hlp}
Suppose that either 
\begin{enumerate}
\item[(a)]
$E[2]$ has a nontrivial Galois-equivariant automorphism  
and $\End_\C(E) \ne \Z[i]$, or
\item[(b)]
$E$ has a rational subgroup of odd prime order $p$ 
and $\End_\C(E) \not\supset \Z[\sqrt{-p}]$.
\end{enumerate}
Then there is a squarefree
polynomial $g(u)$ of degree $6$ such that the 
twist $E_g$ has rank two over $\Q(u)$.
\end{thm}

\begin{proof}
Suppose first that we are in case (a).  
Let $h(t)$ be the linear fractional transformation 
which (after identifying the roots of $f(x)$ with the nonzero 
elements of $E[2]$) agrees with the given automorphism of $E[2]$
on the roots of $f$.  It follows from the Galois-equivariance of 
the automorphism that $h \in \Q(t)$.  
If $h(t) = \alpha t + \beta$, then (since $h(t) \ne t$) 
we must have $\alpha = -1$, and then the set of roots of $f$ 
must be of the form 
$\{\frac{\beta}{2}-a,\frac{\beta}{2},\frac{\beta}{2}+a\}$ for some 
nonzero $a$.  But this contradicts the fact that 
$\End_\C(E) \ne \Z[i]$, so $h$ cannot be a linear polynomial.
Hence in this case the theorem follows from 
Propositions \ref{lfts} and \ref{rk2} and Remark \ref{hlprem}.

Now suppose we are in case (b).  Let $\tilde{E}$ be the quotient of $E$ by 
the given rational subgroup.  Then $\tilde{E}$ is an elliptic curve defined 
over $\Q$ by a 
Weierstrass model $y^2 = \tilde{f}(x)$, and there is an  
isogeny $\phi : \tilde{E} \to E$ of degree $p$, also defined over $\Q$.
Let $h(t) = \phi_x(\mu(t))$ where $\phi_x$ is the 
$x$-coordinate of the isogeny $\phi$ (as in Remark \ref{isogsrem}) 
and $\mu$ is the the linear fractional transformation which 
maps the roots of $f$ to the roots of $\tilde{f}$ in the same way as
the dual isogeny $\hat\phi$ maps $E[2]$ to $\tilde{E}[2]$.  Since 
$\hat\phi$ is defined over $\Q$, $\mu \in \Q(t)$.  
If $\mu(t) = \alpha t + \beta$, then after replacing $\tilde{f}(x)$ by 
$\tilde{f}(x+\beta)$ we may assume that $\beta = 0$.  Then multiplication 
by $\alpha$ sends the roots of $f$ to the roots of $\tilde{f}$, so 
$\tilde{E}$ is the twist of $E$ by $\alpha$.  Let $\iota : E \to \tilde{E}$ 
be an isomorphism over $\C$.  Then $\phi\circ\iota \in \End_\C(E)$ 
and $(\phi\circ\iota)^2 = -p$.  This is impossible since we assumed that
$\sqrt{-p} \notin \End_\C(E)$, so $\mu$ cannot be a linear polynomial.
Now the theorem follows in this case from 
Propositions \ref{isogs} and \ref{rk2} and Remark \ref{hlprem}.
\end{proof}

\begin{rem}
If $E$ has a rational point of order $2$ and $j(E)\ne 1728$,
then hypothesis (a) of Theorem \ref{hlp} holds.
\end{rem}

We illustrate Theorem \ref{hlp} by using the method of 
\S\ref{useful}
to construct some explicit families of examples.
In \S\ref{densect} we will make use of the explicit forms of the
polynomials $g$ below.

If $E$ is an elliptic curve over $\Q$ and
$E(\Q)$ has a point of order $2$, then by translating the $x$-coordinate 
we may assume that $(0,0)$ is a point of order $2$, 
and hence $E$ is of the form $y^2 = x^3 + ax^2 + bx$.

\begin{cor}
\label{rk2deg6}
Suppose that $E$ is $y^2 = x^3 + ax^2 + bx$ with $a, b \in \Q$, 
$ab \ne 0$, $b^2 \ne 4a$.  Let 
$$
g(u) = -ab(u^2 + b^2)(u^4 + 2b^2u^2 - a^2bu^2 + b^4).
$$
Then $E_{g(u)}$ is an elliptic curve over $\Q(u)$ of rank 
$2$, with independent points of infinite order
$$
\Bigl(-\frac{u^2+b^2}{ab},\frac{1}{a^2b^2}\Bigr), \quad
    \Bigl(-\frac{b(u^2+b^2)}{a u^2},\frac{b}{a^2 u^3}\Bigr).
$$
\end{cor}

\begin{proof}
That these points belong to $E_g(\Q(u))$ can be checked directly.  
Since they are nonconstant, both points have infinite order.  
The automorphism of $\Q(u)$ which sends $u$ to $-u$ fixes  
the first point and sends the second point to its inverse, so they 
are independent in $E_{g}(\Q(u))$.  
Since $\deg(g) = 6$, Remark \ref{hlprem} and Lemma \ref{hypgenus} 
show that the rank cannot be greater than two.
\end{proof}

\begin{rem}
Corollary \ref{rk2deg6} was obtained by the method of 
Propositions \ref{rk2} and \ref{lfts} as follows.  
Let $h(t) = -\frac{bt}{at+b}$, the linear fractional transformation 
which switches the two nonzero roots of $f$.  
(This is where we use that $f$ has a rational root; if not,  
$h$ would not have rational coefficients.)
By Propositions \ref{rk2} and \ref{lfts} we see that $E_{f(t)}$ 
has rank at least 2 over $\Q(t,\sqrt{-b(at+b)}) = \Q(u)$ 
where we can take $u = \sqrt{-b(at+b)}$.
We then have $t = -(u^2+b^2)/(ab)$, 
and writing the curve $E_{f(t)}$ and the points 
$(t,1)$, $(h(t),\sqrt{f(h(t))/f(t)})$ in terms of $u$ 
we obtain the data in Corollary \ref{rk2deg6}.
\end{rem}

Suppose now that $E$ has a $\Q$-rational subgroup of order $3$.  The 
$x$-coordinate of the two nonzero points in this subgroup is rational, 
and after translating we may assume that this $x$-coordinate is zero.  
With this normalization one computes that $E$ has a model of the form
$$
y^2 = x^3 + (b^2/4c)x^2 + bx + c
$$
with $b, c \in \Q$, $c \ne 0$, $b^3 \ne 54c^2$, 
and conversely every curve 
defined by such an equation has a $\Q$-rational subgroup 
$\{O,(0,\sqrt{c}),(0,-\sqrt{c})\}$.

\begin{cor}
\label{rk2deg6isog}
Suppose that $E$ is $y^2 = x^3 + (b^2/4c)x^2 + bx + c$ with 
$b, c \in \Q$, $bc \ne 0$, $b^3 \ne 54c^2$.  Let 
$$
g(u) = -bc(2u^6 + (18c^2-b^3)u^4 + (54c^4+2b^3c^2)u^2 + 54c^6-b^3 c^4).
$$
Then $E_{g(u)}$ is an elliptic curve over $\Q(u)$ of rank $2$, 
with independent points of infinite order
$$
\Bigl(-\frac{u^2+3c^2}{2bc},\frac{1}{4b^2c^2}\Bigr), \quad
    \Bigl(\frac{c g(u) - b^4u^2(u^2-c^2)^2)}{4 b^2 c u^2 (u^2+3c^2)^2},
        \frac{c g(u) + 3b^4u^2(u^2-c^2)^2)}{8 b^3 c u^3 (u^2+3c^2)^3}\Bigr).
$$
\end{cor}

\begin{proof}
As with Corollary \ref{rk2deg6}, the simplest proof is a direct 
calculation.
\end{proof}

\begin{rem}
Corollary \ref{rk2deg6isog} was obtained by the method of 
Propositions \ref{rk2} and \ref{isogs} as follows.  
The quotient of $E$ by the subgroup of order $3$ generated 
by $(0,\sqrt{c})$ is the curve $\tilde{E}$ given by $Y^2 = \tilde{f}(X)$ where
$$
\tilde{f}(X) = X^3 - \frac{3b^2}{4c} X^2 - \frac{b(b^3-54c^2)}{6c^2} X 
    - \frac{(b^3-54c^2)^2}{108c^3}.
$$
Let $\phi : \tilde{E} \to E$ be the isogeny 
given by $(\phi_x(X),Y\phi_y(X))$ where
\begin{align*}
\phi_x &= \frac{- 27c^3x^3 + 27b^2c^2x^2 - (9b^4 c - 486bc^3) x
    + b^6 - 108b^3c^2 + 2916c^4}{243c^3x^2}, \\
\phi_y &= \frac{27c^3x^3 - (9b^4 c - 486bc^3) x 
    + 2b^6 - 216b^3c^2 + 5832c^4}{729c^3x^3}.
\end{align*}
The linear fractional transformation $\mu(t)$ which sends the 
roots of $f$ to the roots of $\tilde{f}$ in the same way that $\hat\phi$ 
sends $E[2]$ to $\tilde{E}[2]$ is 
$$
\mu(t) = \frac{(b^3-54c^2)t}{6c(2bt+3c)}.
$$
As in Proposition \ref{isogs} we take $h(t) = \phi_x(\mu(t))$ and 
see that $E_{f(t)}$ has rank two over 
$\Q(t,\sqrt{-c(2b t+3c)}) = \Q(u)$ where we let $u = \sqrt{-c(2b t+3c)}$.
Then $$t = -(u^2+3c^2)/(2bc),$$ 
and writing the curve $E_{f(t)}$ and the points 
$(t,1)$, $(h(t),\sqrt{f(h(t))/f(t)})$ in terms of $u$ 
we obtain the data of Corollary \ref{rk2deg6isog}.
\end{rem}

The following example is contained in Th\'eor\`eme 1 of Mestre \cite{Mestre}.  
We include it here to show how it fits into the framework of this paper.
This result includes the families in Corollaries \ref{rk2deg6}
and \ref{rk2deg6isog} above. The advantages of those corollaries is that the
polynomials $g(u)$ have smaller degree, which will lead to stronger results 
in \S\ref{densect}.

\begin{thm}[\cite{Mestre}]
Suppose that $E : y^2 = x^3 + ax + b$ is an elliptic curve over $\Q$
with $ab \ne 0$. Let 
$$
g(u) = -ab(b^2 (u^4+u^2+1)^3 + a^3 u^4 (u^2+1)^2)(u^2+1).
$$
Then $E_{g(u)}$ has rank at least $2$ over $\Q(u)$.
\end{thm}

\begin{proof}
Let $f(x) = x^3 + ax + b$,
$h_1(t) = -\frac{b(t^3-1)}{a(t^2-1)}$, and 
$h_2(t) = -\frac{b(t^3-1)}{at(t^2-1)}$, and apply Corollary \ref{method} with 
$g_i=f\circ h_i$ and $u = \sqrt{t}$.
\end{proof}

\section{Rank 3}
\label{rank3}

Suppose for this section that $E(\Q)$ contains $3$ points of order $2$, 
i.e., $f(x)$ has three rational roots.  After translating 
and scaling (scaling corresponds to taking a quadratic twist, which is 
harmless for our purposes) we may assume that $f(x) = x(x-1)(x-\lambda)$
with $\lambda \in \Q - \{0, 1\}$.

Suppose $\sigma$ is a permutation of the roots $\{0,1,\lambda\}$ of $f$.
There is a unique linear fractional transformation $h_\sigma(t) \in \Q(t)$ 
which acts on $\{0,1,\lambda\}$ as $\sigma$ does.  
By Proposition \ref{lfts}, 
as long as $h_\sigma(t)$ is not linear there are $j_\sigma \in \Q(t)$ 
and $k_\sigma \in \Q[t]$ 
such that $f\circ h_\sigma = k_\sigma f j_\sigma^2$. 

In order to use these $h_\sigma$ in Proposition \ref{rk3}, we will need 
to find $\sigma_1, \sigma_2$ such that the curve defined by 
$r^2 = k_{\sigma_1}(t)$, $s^2 = k_{\sigma_2}(t)$ has a rational point.

\begin{thm}
\label{2k2}
Suppose that $E$ is an elliptic curve of the form
$y^2=x(x-1)(x-\lambda)$ where
$\lambda = -2a^2$ with $a \in \Q^\times$.
Let $g(u)$ be the polynomial of degree $12$ in $u$
$$
g(u) = 2 N(\lambda,u)(N(\lambda,u) - 2 D(\lambda,u)^2)
    (N(\lambda,u)-2\lambda D(\lambda,u)^2)
$$ 
where 
\begin{align*}
D(\lambda,u) &= \lambda(2\lambda-1)u^2 + 2-\lambda, \\
N(\lambda,u) &=  \lambda^2(\lambda+1)(2\lambda-1)^2 u^4 
   - 4\lambda^2(\lambda-1)(2\lambda-1) u^3 \\
   &\qquad + 2\lambda(\lambda+1)(2\lambda^2-3\lambda+2) u^2 
   - 4\lambda(\lambda-1)(\lambda-2) u
      + (\lambda-2)^2(\lambda+1).
\end{align*}
Then $E_{g(u)}$ has rank at least $3$ over $\Q(u)$, with independent 
points
\begin{align*}
P_1 &= \bigl(\textstyle\frac{N(\lambda,u)}{2 D(\lambda,u)^2},
    \textstyle\frac{1}{4 D(\lambda,u)^3}\bigr), \\
P_2 &= \bigl(\textstyle\frac{\lambda^2 (D(\lambda,u)^2 
    - 4\lambda u (u-1) (\lambda(2\lambda-1)u + 2-\lambda))}
        {(\lambda(2\lambda-1)u^2 - 2\lambda(2\lambda-1)u+\lambda-2)^2},
    \textstyle\frac{a\lambda}
        {(\lambda(2\lambda-1)u^2 - 2\lambda(2\lambda-1)u+\lambda-2)^3}\bigr), \\
P_3 &= \bigl(\textstyle\frac{D(\lambda,u)^2 
    + 4\lambda u (u-1) (\lambda(2\lambda-1)u + 2-\lambda)}
        {\lambda(\lambda(2\lambda-1)u^2 - (2\lambda-4)u 
        + \lambda-2)^2},
    -\textstyle\frac{a}
        {\lambda^2(\lambda(2\lambda-1)u^2 - (2\lambda-4)u + \lambda-2)^3}
    \bigr).
\end{align*}
\end{thm}

\begin{proof}
Take $\sigma_1$ to be the permutation of $\{0,1,\lambda\}$ 
which switches $0$ and $1$, and 
$\sigma_2$ to be the permutation which switches $0$ and $\lambda$.
Then the linear fractional transformations
$$
h_1(t) = \frac{\lambda^2 t - \lambda^2}{(2\lambda-1)t-\lambda^2}, \qquad
    h_{2}(t) = \frac{-t+\lambda}{(\lambda-2)t+1}
$$
act on $\{0,1,\lambda\}$ as $\sigma_1$ and $\sigma_2$ do, respectively. 
One computes in Propositions \ref{lfts} that 
$f\circ h_1 = k_1 f j_1^2$ and 
$f\circ h_{2} = k_{2} f j_{2}^2$ where
$$
k_1(t) = (1-\lambda)((\lambda-2)t+1), \quad 
    k_{2}(t) = \lambda(1-\lambda)((2\lambda-1)t - \lambda^2).
$$
If $a\ne 0$, then $k_{1}$ and $k_{2}$ are $\Q$-linearly independent.
Setting $t_0 = (\lambda+1)/2$, and using  
that $\lambda = -2a^2$, one obtains 
$$
k_1(t_0) = k_{2}(t_0) = a^2 (\lambda-1)^2.
$$
These formulas give us a rational point on the curve of genus zero 
defined by $r^2 = k_1(t), s^2 = k_{2}(t)$.
Using this point one computes that 
$\Q(t,\sqrt{k_1(t)},\sqrt{k_{2}(t)}) = \Q(u)$
where
$$
u = \frac{\sqrt{k_{2}(t)} - a(\lambda-1)}
    {\sqrt{k_{1}(t)} - a(\lambda-1)},
$$
and then $t = N(\lambda,u)/(2 D(\lambda,u)^2)$.  Hence if $g(u)$ is as in 
the statement of the theorem, then $f(t) = g(u)/(4 D(\lambda,u)^3)^2$ 
and the theorem follows from Proposition \ref{rk3}.  
The $3$ points of infinite 
order are computed by taking points
with $x$-coordinates $t$, $h_1(t)$, and $h_{2}(t)$, 
and expressing $t$ in terms of $u$.
\end{proof}

\begin{thm}
\label{1k2}
Suppose that $E$ is given by $y^2=x(x-1)(x-\lambda)$ where either 
\begin{enumerate}
\item[(a)]
$\lambda = \frac{1-a^2}{a^2+2}$ with $a \in \Q-\{0, \pm 1\}$, or
\item[(b)]
$\lambda = \frac{a(a-2)}{a^2+1}$ with $a \in \Q-\{0, 2\}$.
\end{enumerate}
Then there is a squarefree polynomial $g(u) \in \Q[u]$ of degree $12$ in $u$, which 
factors into a product of three quartic polynomials, such that 
$E_{g(u)}$ has rank at least $3$ over $\Q(u)$.
(See \cite{ElecApp} for the polynomials $g(u)$ and 
independent points of infinite order.)
\end{thm}

\begin{proof}
Take $\sigma_1$ to be the permutation of $\{0,1,\lambda\}$ 
which switches $0$ and $1$, $\sigma_2$ to be the permutation 
which switches $1$ and $\lambda$, and 
$\sigma_3$ to be the cyclic permutation $0 \mapsto \lambda \mapsto 1 \mapsto 0$.
Let $h_i \in \Q(t)$ be the corresponding linear fractional 
transformation. Then in Propositions \ref{lfts} we have 
$f\circ h_i = k_i f j_i^2$ where
$$
k_1(t) = (1-\lambda)((\lambda-2)t+1), \, 
    k_2(t) = (1-\lambda)\lambda((\lambda^2-\lambda+1)t - \lambda), \, 
    k_3(t) = \lambda((\lambda+1)t-\lambda).
$$

Now suppose $\lambda = \frac{1-a^2}{a^2+2}$ with $a \in \Q-\{0, \pm 1\}$.
Then $k_1$ and $k_2$ are $\Q$-linearly independent, and
setting $t_0 = \frac{2\lambda}{\lambda+1}$ we find
$$
k_1(t_0) = a^2 (\lambda-1)^2, \quad
    k_2(t_0) = a^2 \lambda^2 (\lambda-1)^2.
$$
These formulas give us a rational point on the curve 
$r^2 = k_1(t), s^2 = k_2(t)$.

If $\lambda = \frac{a(a-2)}{a^2+1}$ with $a \in \Q-\{0, 2\}$,
then $k_2$ and $k_3$ are $\Q$-linearly independent, and 
setting $t_0 = \frac{1}{\lambda}$ we find
$$
k_2(t_0) = (\lambda-1)^2, \quad
    k_3(t_0) = \bigl(\textstyle\frac{a^2+a-1}{a^2+1}\bigr)^2.
$$
These formulas give us a rational point on the curve 
$r^2 = k_2(t), s^2 = k_3(t)$.

The theorem now follows from Proposition \ref{rk3}.
\end{proof}

The following example applies to essentially the same curves 
as Th\'eor\`eme 2 of \cite{Mestre98}.

\begin{thm}
\label{2by4}
Suppose $E[2] \subseteq E(\Q)$ 
and $E$ has a rational cyclic subgroup of order $4$.  
Then $E$ has a model $$y^2 = x(x-b)(x-a^2b)$$ where 
$a,b \in \Q^\times$, $a \ne 1$.
Let $g(u)$ be the poloynomial of degree $11$
\begin{align*}
g(u) &= -4bu\Bigl((a-1)^2u-a\Bigr)
  \Bigl(a^2(a^2-3a+4)u - (a^2+1)(a-1)\Bigr)\\
  &\times\Bigl(a(a^2 - 3a + 4)u^2 -2a(a-1)u + a+1\Bigr)\\
  &\times\Bigl(a(a+1)(a-1)^2(a^2-3a+4)u^2 - 2a(a-1)^2(a^2+1)u + (a^2+1)^2\Bigr)\\
    &\times\Bigl(a^2(a-1)^2(a^2-3a+4)^2u^4 - 4a^2(a-1)^3(a^2-3a+4)u^3   \\ 
 &\qquad + 2(a-1)^2(3a^4-6a^3+5a^2+2)u^2 - 4a(a-1)^2(a^2+1)u + (a^2+1)^2\Bigr).
\end{align*}
Then $E_{g(u)}$ has rank at least $3$ over $\Q(u)$.
(See \cite{ElecApp} for $3$ independent points of infinite order.)
\end{thm}

\begin{proof}
We may write $E$ as $y^2 = f(x)$ where $f$ has $3$ rational roots.  
If $C_4$ denotes the rational cyclic subgroup of order $4$, 
then $2C_4$ contains a rational point, and we may choose our 
model so that this point is $(0,0)$.  Denote the other 
roots of $f$ by $b$ and $b \lambda$.
If $Q$ is a generator of $C_4$ and $x(Q)$ is its $x$-coordinate, 
then $x(Q) \in \Q$ and a computation gives $x(Q)^2 = \lambda b^2$.
Hence $\lambda$ is a square, and we write 
$\lambda = a^2$ with $a \in \Q^\times$.  Thus $E$ is given by 
$y^2 = f(x) := x(x-b)(x-a^2b)$.

The quotient of $E$ by the group generated by $(0,0)$ is 
$$
\tilde{E} : Y^2 = \tilde{f}(X) := X(X+(a-1)^2b)(X+(a+1)^2b).
$$
The isogeny from $\tilde{E}$ to $E$ is 
$\phi(X,Y) = (\phi_x(X),Y\phi_y(X))$ where 
$$
\phi_x(X) = \frac{(X+(a-1)^2b)(X+(a+1)^2b)}{4X}, \quad
    \phi_y(X) = \frac{X^2-(a^2-1)^2b^2}{8X^2}.
$$
The linear fractional transformation 
$$
\mu(t) = \frac{a(a+1)(a-1)^2b (t-b)}{-(a^2-3a+4)t + a(a+1)b}
$$
sends the roots of $f$ to the roots of $\tilde{f}$.
Set $h_1(t) = \phi_x(\mu(t)) \in \Q(t)$.

Let $\sigma$ be the permutation of $\{0,b,a^2b\}$ which switches 
$b$ and $a^2b$, and let $h_2 \in \Q(t)$ be the 
corresponding linear fractional transformation.
One computes in Propositions \ref{isogs} and \ref{lfts} that 
$f\circ h_1 = k_1 f j_1^2$ and $f\circ h_2 = k_2 f j_2^2$ where 
$$
k_1(t) = (a-1)ab((a^2-3a+4)t - a(a+1)b), \quad 
    k_{2}(t) = b((a^2+1)t - a^2b).
$$
Setting $t_0 = a^2b$ we find
$$
k_1(t_0) = (a-1)^4a^2b^2, \quad
    k_{2}(t_0) = a^4b^2.
$$
These formulas give us a rational point on the curve defined by 
$r^2 = k_1(t), s^2 = k_2(t)$.  Using this point one computes that 
$\Q(t,\sqrt{k(t)},\sqrt{k_{\sigma}(t)}) = \Q(u)$
where
$$
u = \frac{\sqrt{k_1(t)} - (a-1)^2ab}
    {\sqrt{k_{2}(t)} - a^2b}.
$$
We can solve for $t$ in terms of $u$ (see \cite{ElecApp}).
The theorem then follows from Proposition \ref{rk3}.
\end{proof}

\begin{rem}
The theorems above give certain infinite families of curves which 
have twists of rank (at least) 3 over $\Q(u)$.  The restriction to 
these families makes it possible to find rational points on the 
genus zero curves $r^2 = k_1(t), s^2 = k_2(t)$ which arise in the 
construction.  It is possible to carry out the construction for many 
curves not in these families.  We give one example in the 
next theorem.
\end{rem}

\begin{thm}
\label{cnum}
The elliptic curve
$$
6(u^{12}-33u^8-33u^4+1)y^2=x^3-x
$$ 
has rank at least $3$ over $\Q(u)$, with independent points
$$
P_1 = \bigl(-\textstyle\frac{u^4 - 6u^2 + 1}{3(u^2+1)^2}, 
    \textstyle\frac{2}{9(u^2+1)^3}\bigr), \quad
P_2 = \bigl(-\textstyle\frac{u^4 + 6u^2 + 1}{3(u^2-1)^2}, 
    \textstyle\frac{2}{9(u^2-1)^3}\bigr), \quad
P_3 = \bigl(\textstyle\frac{u^4 + 1}{6u^2}, 
    \textstyle\frac{1}{36u^3}\bigr).
$$
\end{thm}

\begin{proof}
The simplest proof is a direct computation.  To construct this 
example one takes $E$ to be $y^2 = x^3 - x$ and 
proceeds exactly as in the proofs of Theorems \ref{2k2} and \ref{1k2}, 
with $h_1(t) = \frac{t-1}{3t+1}$, $h_2(t) = \frac{t+1}{3t-1}$, 
which gives
$$
k_1(t) = 6t+2, \quad k_2(t) = -6t+2.
$$
The curve defined by $r^2 = k_1(t), s^2 = k_2(t)$ has a rational point 
$(r,s,t) = (2,0,1/3)$, and using this one computes that 
$\Q(t,\sqrt{k_1},\sqrt{k_2}) = \Q(u)$ where $t = u^4-6u^2+1$.  
Proposition \ref{rk3} with this input leads to the data above.
\end{proof}

\begin{rem}
\label{favoritecurve}
Let $g(u)=6(u^3-33u^2-33u+1)$. Over $\Q(u)$, the rank of
$E_g$ (respectively $E_{g(u^2)}$, respectively $E_{g(u^4)}$)
is $1$ (respectively $2$, respectively $3$).
Unfortunately this pattern does not continue; the rank of
$E_{g(u^8)}$ is $3$. 
Replacing $u$ by $\sqrt{u}$ in $P_1$ and $P_2$ above
gives two independent points on $E_{g(u^2)}$.
\end{rem}

\section{Densities}
\label{densect}

Recall the definitions of $r_E(D)$ and  $N_r(E,x) \ge N_{r}^+(E,x)$
from the introduction.
In this section 
we use results of Stewart and Top \cite{STop} to obtain lower bounds 
for $N_r(E,x)$ (and, subject to the Parity Conjecture, for
$N_{r}^+(E,x)$, as Gouv\^ea and Mazur \cite{Gouvea-Mazur} did), 
with $E$ and $r$ provided by the examples of the 
previous sections.  The first two assertions of the following theorem 
are immediate from Theorems 2 and 1 of \cite{STop}, and were used by Stewart 
and Top in that paper in several families of examples.  What is new 
here is that by using the examples of the previous sections we have 
more curves to which we can apply these results.  In addition, we 
show in Theorem \ref{gd}(iii) 
how to use Theorem 1 of \cite{STop} 
along with the Parity Conjecture to obtain results for higher rank 
(see also \cite{Gouvea-Mazur} and \S12 of \cite{STop}).

If $A$ is an elliptic curve over $\Q$, let $w(A) \in \{\pm 1\}$ 
denote the root number 
in the functional equation of the $L$-function $L(A,s)$.  The 
Parity Conjecture asserts that $w(A) = (-1)^{\rank(A(\Q))}$.

\begin{thm}
\label{gd}
Suppose that $E$ is an elliptic curve over $\Q$, 
and $g \in \Q[u]$ is nonconstant and squarefree. 
Let $r = \rank(E_g(\Q(u)))$ and $k = \bigl[\frac{\deg(g)+1}{2}\bigr]$.
\begin{enumerate}
\item
For $x\gg 1$,
$$N_r(E,x) \gg x^{1/k}/\log^2(x).$$
\end{enumerate}
Suppose further that the irreducible factors of $g$ all have degree at most $6$.
\begin{enumerate}
\setcounter{enumi}{1}
\item
For $x\gg 1$,
$$N_r(E,x) \gg x^{1/k}.$$
\item
Suppose that the Parity Conjecture holds for all twists of $E$, 
and that there is a rational 
number $c$ such that $g(c) \ne 0$ and  
$w(E_{g(c)}) = (-1)^{r+1}$.  
Then for $x\gg 1$,
$$N_{r+1}^+(E,x) \gg x^{1/k}.$$
\end{enumerate}
\end{thm}

\begin{proof}
Without loss of generality we may assume that
$\deg(g)\ge 3$, since if not, $r=0$ by
Remark \ref{hlprem} and there is nothing to prove.

Let $F(X,Y) = Y^{2k}g(X/Y)$, a homogeneous polynomial 
of degree $2k$. 
Assertions (i) and (ii) are immediate from Theorems 2 and 1 of \cite{STop},
respectively, applied to $F$. 

Suppose now that the Parity Conjecture holds,
the irreducible factors of $g$ all have degree at most $6$,
and $c \in \Q$ 
is such that $g(c) \ne 0$ and $w(E_{g(c)}) = (-1)^{r+1}$.  
Choose a closed interval $I \subset \R$ with rational endpoints which 
contains $c$ but does not contain any roots of $g$, and 
let $\mu(u) = \frac{\alpha u + \beta}{\gamma u + \delta} \in \Q(u)$ 
be a linear fractional transformation which maps $[0, \infty]$ onto $I$
and (for simplicity) such that $\mu(1) = c$.  Replace $g$ by 
the polynomial $(\gamma u + \delta)^{2k}(g \circ \mu)$ of degree 
at most $2k$.  Then we still have that $r = \rank(E_g(\Q(u)))$, 
and our construction guarantees that this new 
polynomial $g$ also satisfies:
\begin{enumerate}
\item[(a)]
the constant term of $g$ and the coefficient of $u^{2k}$ are both nonzero,
\item[(b)]
the irreducible factors of $g$ have degree at most $6$,
\item[(c)]
$g(1) \ne 0$ and $w(E_{g(1)}) = (-1)^{r+1}$,
\item[(d)]
$g(u)/g(1)$ is positive if $u \ge 0$.
\end{enumerate}
Further, multiply $g$ by the square of an integer to
clear denominators of the coefficients.
If $A$ is an elliptic curve over $\Q$, write
$\cond(A)$ for its conductor. If further 
$D \in \Q^\times$ and
$\cond(A)$ is relatively prime to 
the conductor of the character $\chi_D$ associated
to the quadratic extension $\Q(\sqrt{D})/\Q$, then
$w(A_D) = \chi_D(-\cond(A))w(A)$. Applying this with 
$A = E_{g(1)}$ and $D = g(a/b)/g(1)$ for $a$ and $b$
positive integers congruent to $1$ modulo an integer $M$
sufficiently divisible by the prime divisors of $2\cond(E_{g(1)})$, 
and using (c) and (d) above, gives that 
\begin{equation}
\label{sign}
w(E_{g(a/b)}) = w(E_{g(1)}) = (-1)^{r+1}.
\end{equation}
Let 
\begin{multline*}
S = \{\text{squarefree integers $D$} : 
    \text{$D = F(a,b)/v^2$ for some $a,b,v \in \Z^+$} \\
    \text{with
    $a,b \le x$, $a \equiv b \equiv 1 \pmod{M}$}\}, 
\end{multline*}
$$
S(x) = \{D \in S : |D|<x\}.
$$
By Theorem 1 of \cite{STop},  
for $x\gg 1$,
\begin{equation}
\label{st}
\#(S(x)) \gg x^{1/k}.
\end{equation}
(Note that as stated, Theorem 1 of \cite{STop} does not include the 
restriction $a, b > 0$ in our definition of $S(x)$.  However, 
the proof in \cite{STop} does restrict to positive $a, b$.)

Theorem C of \cite{Silvermanhts} implies that $r_E(D) \ge r$ for all 
but finitely many $D \in S$.
However, by \eqref{sign}, if $D \in S$ then $w(E_D) = (-1)^{r+1}$ 
so the Parity Conjecture tells us that $r_E(D) \ne r$.  
Hence $r_E(D) \ge r+1$ for all 
but finitely many $D \in S$, and so assertion (iii) of the theorem 
follows from the Stewart-Top bound \eqref{st}.
\end{proof}

\begin{cor}
\label{rohrcor}
Suppose that $E$ is an elliptic curve over $\Q$, 
and $g \in \Q[u]$ is a nonconstant squarefree polynomial 
whose irreducible factors have degree at most $6$.  
Let $r = \rank(E_g(\Q(u)))$ and $k = \bigl[\frac{\deg(g)+1}{2}\bigr]$.
If the Parity Conjecture holds for all 
twists of $E$, and $g$ has at least 
one real root, then for $x\gg 1$,
$$
N_{r+1}^+(E,x) \gg x^{1/k}.
$$
\end{cor}

\begin{proof}
If $g$ has a real root then $g(\Q)$ contains both positive and negative 
values ($g$ has no multiple roots because it was 
assumed to be squarefree).  
Thus by a result of Rohrlich (Theorem 2 of \cite{Rohrlich}) 
we have 
$$
\{w(E_{g(a)}) : a \in \Q, g(a) \ne 0\} = \{1, -1\}.
$$
Now the corollary follows immediately from Theorem \ref{gd}(iii).
\end{proof}

We now give some applications of Theorem \ref{gd} and 
Corollary \ref{rohrcor}.

\begin{thm}
\label{hlpcor}
Suppose that either 
\begin{enumerate}
\item[(a)]
$E[2]$ has a nontrivial Galois-equivariant automorphism 
and $\End_\C(E) \ne \Z[i]$, or
\item[(b)]
$E$ has a rational subgroup of odd prime order $p$ 
and $\End_\C(E) \not\supset \Z[\sqrt{-p}]$.
\end{enumerate}
Then for $x\gg 1$,
$$
N_2(E,x) \gg x^{1/3}.
$$
\end{thm}

\begin{proof}
This is immediate from Theorems \ref{hlp} and \ref{gd}(ii).
\end{proof}

\begin{thm}
\label{hlpcor2}
Suppose that 
$E[2] \subset E(\R)$ and either 
\begin{enumerate}
\item[(a)]
the largest or smallest root of $f$ is rational, or
\item[(b)]
$E$ has a rational subgroup of order $3$.
\end{enumerate}
If the Parity Conjecture holds for all twists of $E$ then 
for $x\gg 1$,
$$
N_{3}^+(E,x) \gg x^{1/3}.
$$
\end{thm}

\begin{proof}
Suppose first that we are in case (a).  Translating the 
given rational root 
of $f$ we may assume that $f(x) = x^3 + ax^2 + bx$ with $b > 0$. Since 
$f$ has $3$ real roots we also have $a^2 - 2b > a^2 - 4b > 0$.  
In particular, $b(a^2 - 2b) > 0$.
Let $g(u)$ be as in Corollary \ref{rk2deg6}.  Then $g$ is divisible by 
$g_1(u) = u^4 - b(a^2-2b)u^2 + b^4$.  We compute that 
$g_1(\sqrt{b(a^2 - 2b)}) = -\frac{1}{4}a^2b^2(a^2-4b) < 0$, but 
$g_1(u)$ is positive for large $u$, so $g_1$, and hence $g$, has 
real roots.  Hence the Corollary in this case follows from 
Corollaries \ref{rohrcor} and \ref{rk2deg6}.

Similarly, suppose we are in case (b).  Then as discussed before 
Corollary \ref{rk2deg6isog}, $E$ has a model 
$y^2 = x^3 + (b^2/4c)x^2 + bx + c$ with $b, c \in \Q$, $c \ne 0$.  
The discriminant of this model is $\Delta(E) = 8(b^3 - 54c^2)$.  Since 
all the $2$-torsion on $E$ is defined over $\R$, we have $\Delta(E) > 0$.
Let $g(u)$ be as in Corollary \ref{rk2deg6isog}.  Then 
$g(u)/(-bc)$ is positive for large $u$, but 
$g(0)/(-bc) = -c^4(b^3-54c^2) = -c^4\Delta(E)/8 < 0$.
Hence $g$ has real roots, so the Corollary in this case follows from 
Corollaries \ref{rohrcor} and \ref{rk2deg6isog}.
\end{proof}

\begin{thm}
\label{rk3dens}
Suppose $E$ is defined by $y^2 = x(x-1)(x-\lambda)$ where either 
$\lambda = -2a^2$, or
$\lambda = \frac{1-a^2}{a^2+2}$, or  
$\lambda = \frac{a(a-2)}{a^2+1}$,
with $a \in \Q$ and $\lambda \ne 0$. Then for $x\gg 1$,
$$
N_3(x) \gg x^{1/6}.
$$
\end{thm}

\begin{proof}
This is immediate from Theorems \ref{gd}(ii), \ref{2k2}, \ref{1k2}, 
and \ref{cnum} (the last to handle the 
excluded value $a=0$ in Theorem \ref{1k2}(a)).
\end{proof}

\begin{thm}
\label{lastthm}
Suppose $E[2] \subseteq E(\Q)$ and $E$ has a rational cyclic subgroup 
of order $4$.  Then:
\begin{enumerate}
\item
 for $x\gg 1$,
$$
N_3(x) \gg x^{1/6},
$$
\item
if the Parity Conjecture holds for all twists of $E$, then 
for $x\gg 1$,
$$
N_{4}^+(x) \gg x^{1/6}.
$$
\end{enumerate}
\end{thm}

\begin{proof}
Assertion (i) follows directly from Theorems \ref{gd}(ii) and \ref{2by4}.
The polynomial $g$ of Theorem \ref{2by4} has degree $11$, and hence 
it has a real root, so (ii) follows from Corollary \ref{rohrcor} 
and Theorem \ref{2by4}.
\end{proof}

\begin{rem}
The conclusions of Theorems \ref{hlpcor}, \ref{hlpcor2},
\ref{rk3dens}, and \ref{lastthm}
hold when $E$ is $y^2=x^3-x$, by Remark \ref{favoritecurve},  
Theorem \ref{gd}, and Corollary \ref{rohrcor}.
\end{rem}

\section{Remarks and questions}

\theoremstyle{plain}
\newtheorem{prob}[thm]{Problem}
\begin{prob}
\label{prob1}
Find a hyperelliptic curve $C$ of the form $s^2 = g(u)$ with $g(u) \in \Q[u]$ 
such that the jacobian of $C$ is isogenous over $\Q$ to $E^r \times B$ 
for some elliptic curve $E$ and abelian variety $B$, either with $r \ge 4$, 
or with both $r=3$ and $\dim(B)\le 1$.
\end{prob}

\begin{rem}
A solution $(C, E, r, B)$ to Problem \ref{prob1} would imply, by 
Theorem \ref{gd}(i) and the equality in Remark \ref{hlprem}, that 
$$
N_r(E,x) \gg \frac{x^{1/(1+\genus(C))}}{\log^2(x)}
= \frac{x^{1/(1+r+\dim(B))}}{\log^2(x)}.
$$
\end{rem}

\begin{rem}
The reason for the restriction on $r$ in Problem \ref{prob1}
is that we already have examples when $r\le 3$. 
Theorem \ref{hlp} gives numerous examples with $r=2$ and $\dim(B)=0$, and
Theorems \ref{2k2}, \ref{1k2}, \ref{2by4}, and \ref{cnum} 
provide numerous examples with $r=3$ and $\dim(B)=2$.
\end{rem}

\begin{rem}
The results of Stewart and Top \cite{STop} would not be
needed in the arguments of \S\ref{densect} if the following conjecture of 
Caporaso, Harris, and Mazur \cite{CHM} were known to hold.  More precisely,
Proposition \ref{remprop} shows that \eqref{st} above follows 
easily from this conjecture.
\end{rem}

\begin{conj}[Caporaso, Harris, Mazur]
\label{chm}
Fix an integer $h \ge 2$.  Then there is a constant $B(h)$ such that for 
every curve $C$ of genus $h$ defined over $\Q$, $\#(C(\Q)) < B(h)$.
\end{conj}

\begin{prop}
\label{remprop}
Suppose $g(u) \in \Z[u]$ is a squarefree polynomial, and let 
$k = \bigl[\frac{\deg(g)+1}{2}\bigr]$ and $F(X,Y) = Y^{2k}g(X/Y)$.  
Fix a positive integer $M$ and
define $S(x)$ as in the proof of Theorem \ref{gd}(iii), with this $M$.  
If Conjecture \ref{chm} is true and $k \ge 3$, 
then for $x\gg 1$,
$$
\#(S(x)) \gg x^{1/k}.
$$
\end{prop}

\begin{proof}
If $a, b \in \Z$ and $F(a,b) \ne 0$, let $s(F(a,b))$ denote the squarefree part 
of $F(a,b)$,
i.e., the unique squarefree integer $D$ such that $F(a,b) = Dn^2$ for some 
integer $n$.  
For every squarefree integer $D$ let 
$A_D$ denote the hyperelliptic curve 
$Dv^2 = g(u)$ of genus $k-1 \ge 2$.
The map $(a,b) \mapsto (a/b,\pm b^{-k}\sqrt{F(a,b)/D})$ 
defines an injection
$$
\{(a,b) \in \Z^2 : (a,b) = 1, s(F(a,b)) = D\} \hookrightarrow A_D(\Q)/\{\pm 1\}
$$
(where $-1$ denotes the hyperelliptic involution on $A_D$).  
Thus by Conjecture \ref{chm}
the order of the set on the left is bounded by $B(k-1)$.  
Let 
$$
R(x) = \{(a,b) \in \Z^2 : 1 \le a, b \le x, \; (a,b) = 1, \;
F(a,b) \ne 0, \; a \equiv b \equiv 1 \hskip-7pt\pmod{M}\}.
$$
There is a constant $K = K(g)$ such that if $(a,b) \in R(x)$ then $|F(a,b)| < K x^{2k}$.
It follows that $\#(S(x)) \ge \#(R((x/K)^{1/{2k}}))/B(k-1)$ for $x\gg 1$.  But it 
is standard to show that  
$\#(R(x)) \gg x^2/M^2$ for $x\gg 1$, and the 
proposition follows.
\end{proof}

\end{document}